\documentclass[12pt, twoside]{amsart}

\swapnumbers
\numberwithin{equation}{section}

\usepackage{amssymb,latexsym, amsmath, amscd}

\input scrload

\newtheorem{thm}{Theorem}[section]
\newtheorem{theorem}[thm]{Theorem}
\newtheorem{lem}[thm]{Lemma}
\newtheorem{prop}[thm]{Proposition}
\newtheorem{cor}[thm]{Corollary}
\newtheorem{cory}[thm]{Corollary}

\newtheorem{df}[thm]{Definition}

\newcommand\theoref{Theorem~\ref}
\newcommand\lemref{Lemma~\ref}
\newcommand\propref{Proposition~\ref}
\newcommand\corref{Corollary~\ref}

\def \sirc{{\raise0.2ex \hbox{$\scriptstyle \circ$}}}

\def \wt{\widetilde}

\def\ts{\times}

\def\po{\parindent 0pt}
\def\p{{\po \it \s Proof. }}

\def\s{\smallskip}
\def\m{\medskip}

\def\ds{\displaystyle}

\def\cat{\operatorname{cat}}
\def\cd{\operatorname{cd}}
\def\hd{\operatorname{hd}}

\def\C #1{\mathcal #1}

\long\def\forget#1\forgotten{} %

\begin{document}
%\topmargin -2cm

%\vskip 2cm
\title[Cohomological Dimension, Connectivity, and LS category]{Cohomological Dimension, Connectivity, and Lusternik--Schnirelmann category}

%\date{\today}

\author{Yu. B. Rudyak}
\address{Department of Mathematics, 1400 Stadium Rd
University of Florida
Gainesville, FL 32611, USA} 
\email{rudyak@ufl.edu}

\begin{abstract} 
Dranishnikov~\cite{D2} proved that 
\[
 \cat X\leq \cd(\pi_1(X))+\Bigl\lceil\frac{\hd (X)-1}{2}\Bigr\rceil.
\]
where $\cd(\pi)$ denotes the cohomological dimension of a group $\pi$ and $\hd(X)$ denotes the homotopy dimension of $X$. Furthermore, there is a well-known inequality of Grossman,~\cite{G}:
\[
\cat X\leq \Bigl\lceil\frac{\hd (X)}{k+1}\Bigr\rceil \text{ if } \pi_i(X)=0 \text{ for } i\leq k.
\]

\noindent We make a synthesis and generalization of both of these results, by demonstrating the main result:
\[
 \cat X\leq \cd(\pi_1(X))+\Bigl\lceil\frac{\hd (X)-1}{k+1}\Bigr\rceil \text { if  }\pi_i(X)=0 \text{  for } i=2, \ldots, k. 
\]

\noindent The proof of the main theorem uses the Oprea--Strom  inequality
$\cat X\leq \hd(B\pi_1(X))+\cat^1X$,~\cite{OS} where $\cat^1$ is the Clapp-Puppe $\cat \C A$ with $\C A$ the class of 1-dimensional CW complexes. The inequality clarified the Dranishnikov inequality.

\end{abstract}

%\small {\it 2010 Mathematics Subject Classification:}

\maketitle

\section{Introduction}
We work in the category of connected CW complexes and continuous maps. We use the sign $\cong$ for homotopy equivalences. All covers are assumed to be open. Given a space $X$, $\hd(X)$ denotes the {\em homotopy dimension} of $X$, that is, the minimum cellular dimension of all CW complexes homotopy equivalent to $X$. Given a group $\pi$, $B\pi$ denotes a {\em classifying space} for $\pi$, and $\cd(\pi)$ denotes the {\em cohomological dimension} of $\pi$,~\cite{B}.  A {\em classifying map} for $X$ is a map $c=c_X: X \to B\pi_1(X)$ that induces an isomorphism of fundamental groups. 

\m  Below $\cat X$ denotes the{\em  Lusternik--Schnirelmann category} of a space $X$,~\cite{LS, CLOT}. A well-known inequality $\cat X \leq \hd X$,~\cite{LS,F,CLOT} can  be generalized as follows~\cite{G, CLOT}:
\begin{theorem}\label{t:gros}
If
$\pi_i(X)=0$ for $i\leq k$ then $\ds\cat X\leq \frac{\hd X}{k+1}.$
\end{theorem}
Concerning the case of $\pi_1(X)\neq 0$, the author conjectured that $\cat X$ can be asymptotically bounded above by $\hd(X)/2$, provided that $\pi_1(X)$ has finite cohomological dimension, i.e., $\cd(\pi_1(X))<\infty$.
Later  Dranishnikov~\cite{D2} proved the following fact:
\begin{theorem}\label{t:hdim}
$
\ds\cat X\leq \cd(\pi_1(X))+\Bigl\lceil\frac{\hd(X)-1}{2}\Bigr\rceil.
$
\end{theorem}

\m This theorem can be regarded as a confirmation of the conjecture.

\m Also, I suspected that there should be a synthesis of \theoref{t:gros} and \theoref{t:hdim}, so that the equation~\theoref{t:hdim} can be improved by replacing (approximately) $\hd X$ by $\hd X/(k+1$ for $X$ $k$-connected. In other words, I expected to have a claim that generalizes both \theoref{t:gros} and \theoref{t:hdim}. Now I know how to make this improvement (synthesis, generalization). Let me tell you the precise statements (\theoref{main} and \corref{c:cd} below).

\begin{df}[Clapp and Puppe~\cite{CP}]\rm
 Given a class $\mathcal A$ of  CW complexes and a space $X$, define a subset $A$ of $X$ to be {\em $\mathcal A$-categorical} if the inclusion $A\to X$ factors, up to homotopy, through a space in $\mathcal A$. Follow Clapp and Puppe~\cite{CP}, define the $\mathcal A$-cover of $X$ to be the cover $\{U_0, U_1, \ldots, U_m\}$ such that each $U_i$ is $\mathcal A$-categorical. Define  $\cat_{\mathcal A} X$, the {\em $\mathcal A$-category of X } to be the minimal number $k$ such that there exists an  $\mathcal A$-categorical cover $\{U_0, U_1, \ldots, U_k\}$ . 
 \end{df}

\m For example, $\cat_{\mathcal A}(X)=\cat (X)$ if $\mathcal A$ is the class of contractible spaces.

\m \begin{df}\rm Let $\mathcal A(n)$ be the class of all $n$-dimensional CW complexes. Put $\cat^n(X):=\cat_{\mathcal A(n)}(X)$. 
\end{df}

\m The following Oprea--Strom Theorem recovers and clarifies \theoref{t:hdim}.

\begin{theorem}\label{t:os}
For every space $X$ we have
\[
 \cat X\leq \hd (B\pi_1(X))+\cat^1(X) \leq \hd(B\pi_1(X))+\Bigl\lceil\frac{\hd (X)-1}{2}\Bigr\rceil.
 \]
\end{theorem}

 \p See  Oprea and Strom~\cite[Corollary 6.2]{OS}.
 \qed
 
\m We prove the following generalization of \theoref{t:os}:
\begin{thm}[Corollary \ref{c:main}]\label{main}
Let $X$ be a CW complex and $\pi=\pi_1(X)$. Suppose that the classifying map $c: X\to B\pi$ induces an isomorphism $c_*: \pi_i(X)\to \pi_i(B\pi)$ for $i\leq k$. $($In particular, $\pi_i(X)=0$ for $i=2, \ldots, k.)$ Then
\[
 \cat X\leq \hd(B\pi)+\cat^k(X) \leq \hd(B\pi)+\Bigl\lceil\frac{\hd (X)-1}{k+1}\Bigr\rceil.
 \]
\end{thm}

\begin{cor}\label{c:cd}
Let $X$ be a CW complex as in Theorem $\ref{main}$. Then
\[
 \cat X\leq \cd(\pi)+\Bigl\lceil\frac{\hd (X)-1}{k+1}\Bigr\rceil.
 \]
\end{cor}

Clearly, \theoref{main} and \corref{c:cd} can be regarded as an above-mentioned synthesis.

\section{Proofs}

First, we settle the second part of the inequality noted in \theoref{main}. 

\begin{prop}\label{p:second}
Let $X_k$ be the $k$-skeleton of a CW complex $X$. Then
\[
\cat^k(X)\leq  \Bigl\lfloor\frac{\hd (X)}{k+1}\Bigr\rfloor \leq\Bigl\lceil\frac{\hd (X)-1}{k+1}\Bigr\rceil
\]
\end{prop}

\p For the first inequality, see~\cite[Proposition 4.4]{OS}. The second inequality is obvious.
\qed

\m Now we prove the first part of \ref{main}. The proof is based (speculated) on~\cite[Sections 5,6]{OS} that, in turn, exploits clever ideas of Dranishnikov~\cite{D1,D2}. 

\m Let $X$ be a CW complex. Take $k>1$, let $X_k$ be the $k$-skeleton of $X$, and let $\wt Z$ be the universal covering of $Z$ for $Z=X$ or $Z=X_k$. Put $\pi=\pi_1(X)$ and let $E\pi\to B\pi$ be the universal bundle for $\pi$. Note that $\pi$ acts on $\wt Z$ via deck transformations of the covering $\wt Z \to Z$, and we can form the Borel construction 
\begin{equation}\label{eq:borel}
p :E\pi\ts_{\pi}\wt Z \to B\pi.
\end{equation} 
It is worth noting that $E\pi$ is contractible, and so 
\begin{equation}\label{eq:cong}
E\pi\ts_{\pi}\wt Z \cong  Z.
\end{equation}
The inclusion $i=i_k: X_k\to X$ yields the commutative diagram
\[
\CD
E\pi\ts_{\pi}\wt X_k @>f>> E\pi\ts_{\pi}\wt X\\
@Vp_0VV  @VVp_1V\\
B\pi @=B\pi
\endCD
\]
where $f=f_k$ is induced by $i$, and $p_0=p$ if $Z=X_k$, and $p_1=p$ if $Z=X$. Take a base point $*$ of $B\pi$ and let $F_0, F_1$ be the fibers of $p_0, p_1$ over $*$, respectively. Then $f$ yields a map (inclusion) $j: F_0\to F_1$ of fibers. 

\begin{prop}\label{p:conn}
If $\pi_i(X)=0$ for $i=2, \ldots, k$ then the inclusion $j: F_0\to F_1$ is null-homotopic.
\end{prop}

\p It follows because $F_0$ is homotopy equivalent to $\wt X_k$, while $\wt X_k$ is contractible $(\pi_i(\wt X_k)=0$ for $i\leq k$ and $\hd X_k\leq k$).
\qed 

\m Following~\cite{OS}, define a cover $\C U=\{U_0, U_1, \ldots, U_n\}$ of $E\pi\ts_{\pi}\wt X$ to be a $\Gamma_k$-cover if each inclusion $U_m\subset E\pi\ts_{\pi}\wt X$ passes through the inclusion $f$, up to homotopy over $B\pi$. In this case we define $\gamma(\C U)=n$. Now, set
\[
\Gamma_k(X)=\inf\{\gamma(\C U)\bigm | \C U \text{ is a $\Gamma_k$-cover of  } E\pi\ts_{\pi}\wt X\}.
\]

\begin{prop}\label{p:j}
We have $\cat X=\cat (E\pi\ts_{\pi}\wt X)\leq \hd B(\pi)+\Gamma_k$.
\end{prop}

\p The equality $\cat X=\cat (E\pi\ts_{\pi}\wt X)$ is explained in \eqref{eq:cong}. The inequality follows  from~\cite[Prop. 5.1]{OS} because of~\propref{p:conn}.
\qed
 
\begin{thm}\label{t:gamma}
For any CW complex $X$ and every $k\geq 1$, we have $\Gamma_k(X)=\cat^k(X)$
\end{thm}

\p For $k=1$, this is~\cite[Theorem 6.1]{OS}. For $k>1$, the proof is literally the same as for $k=1$. The only change is to replace $X_1$ by $X_k$, $\cat^1$ by $\cat^k$, and $\Gamma_1$ by $\Gamma_k$, in  ~\cite[Theorem 6.1]{OS}. 
\qed

\begin{cory}\label{c:main}
Let $X$ be a CW complex and $\pi=\pi_1(X)$. Suppose that the classifying map $c: X\to B\pi$ induces an isomorphism $c_*: \pi_i(X)\to \pi_i(B\pi)$ for $i\leq k$. Then
\[
\cat X\leq \hd (B\pi_1(X))+\cat^k(X) \leq \hd(B\pi_1(X))+\Bigl\lceil\frac{\dim (X)-1}{k+1}\Bigr\rceil.
\]
\end{cory}

\p The first inequality follows from \propref{p:j} and \theoref{t:gamma}, the second inequality follows from \propref{p:second}.
\qed

\m Now we prove {\bf \corref{c:cd}}. First, given a group $\pi$, recall that $\cd(\pi)=\hd(B\pi)$ if either $\cd(\pi)\leq3$,~\cite{EG} or $\cd(\pi)=1$,~\cite{Stal, Swan}. Furthermore, recall that $\cd(\pi)=\cat(B\pi)$ for all groups $\pi$,~\cite{EG, Stal, Swan}. {\em So, for $\cd(\pi)\neq 2$, \corref{c:cd} follows from \corref{c:main} directly.} 
(Note also that if $\cd(\pi)=2$ then either $\hd (B\pi)=2$ or  $\hd (B\pi)=3$, and it is unknown question whether there exists a group $\pi$ with $\cd(\pi)=2$ and $\hd(B\pi)=3$.)

\m Consider a space $X$ and the classifying map $c: X \to B\pi$ where $\pi=\pi_1(X)$. Note that $c_*: \pi_1(X)\to \pi_1(B\pi)$ is an isomorphism. Given $k\in \mathbb N$, assume that $c_*: \pi_i(X)\to \pi_i(B\pi)$ is an isomorphism for $i\leq k$.

\begin{lem}\label{l:dran}
Let $f: X\to Y$ be a locally trivial bundle with an $k$-connected
fiber $F$. Suppose that $f$ admits a section. Then
\[
\cat X\leq \cat Y+\Bigl\lceil\frac{\hd (X)-k}{k+1}\Bigr\rceil. 
\]
\end{lem}

\p See~\cite[Theorem 3.7]{D1}.
\qed

\m We apply \lemref{l:dran} to the Borel construction $p: E\pi\ts_\pi{\wt X}\to B\pi$ as in \eqref{eq:borel},  with $\cd(\pi)=2$. Note that the bundle $p$ is the classifying map for $X$. Furthermore, the fiber $F$ of $p$ is homotopy equivalent to $\wt X$. 

\m For $k=1$, \corref{c:cd} is the Dranishnikov theorem \theoref{t:hdim}. So, assume that $k>1$. Then $\pi_2(F)=\pi_2(\wt X)=0$, since $\pi_2(\wt X)=\pi_2(X)=\pi_2(B\pi)$. 

\m We have $\hd(B\pi)\leq 3$ and $\pi_i(F)=0$ for $i=1, 2$. So, because of the elementary obstruction theory, $p$ has a section. Thus, because of \ref{l:dran} and since $\cd(\pi)=\cat(B\pi)$, we conclude that 
\[
\cat X \leq \cd \pi +\Bigl\lceil\frac{\hd (X)-1}{k+1}\Bigr\rceil
\]
for $\cd(\pi)=2$, and therefore for all $\pi$. This completes the proof of \corref{c:main}.

\m {\bf Acknowledgments:} The work was partially supported by a grant from the Simons Foundation (\#209424 to Yuli Rudyak).

\end{document}